\documentclass[notes-porbanz,noinfoline]{notes-imsart}
\pdfoutput=1

\RequirePackage[OT1]{fontenc}
\RequirePackage{amsthm,amsmath}
\usepackage[usenames,dvipsnames]{color}
\definecolor{RefColor}{rgb}{0,0,.65}
\RequirePackage[colorlinks,linkcolor=RefColor,citecolor=RefColor,urlcolor=RefColor]{hyperref}

\setattribute{copyright}{text}{\empty}
\hypersetup{pdfinfo={Subject={ }}}

\usepackage{amsmath,amssymb,amscd,amsfonts,amsthm}
\usepackage{thmtools}
\usepackage[numbers]{natbib}
\usepackage[capitalize]{cleveref}
\usepackage{booktabs}
\usepackage{tikz}
\usetikzlibrary{shapes.geometric,calc,arrows,chains,matrix,positioning,scopes,decorations.markings,decorations.pathreplacing,intersections,backgrounds,fit,scopes}
\usepackage{dsfont}
\usepackage{mathtools}

\startlocaldefs
\theoremstyle{plain}
\newtheorem{theorem}{Theorem}
\newtheorem{lemma}[theorem]{Lemma}
\newtheorem{proposition}[theorem]{Proposition}

\newtheorem{corollary*}[theorem]{Corollary}

\theoremstyle{definition}

\makeatletter
\newcommand{\eqnum}{\leavevmode\hfill\refstepcounter{equation}\textup{\tagform@{\theequation}}}
\makeatother
\endlocaldefs

\begin{document}
\begin{abstract}
  We give conditions under which a scalar random variable $T$ can be coupled to a random
  scaling factor $\xi$ such that $T$ and ${\xi T}$ are rendered stochastically independent.
  A similar result is obtained for random measures. 
  One consequence is a generalization of a result by Pitman and Yor
  on the Poisson-Dirichlet distribution to its negative parameter range.
  Another application are diffusion excursions straddling an exponential random time.
\end{abstract}

\begin{frontmatter}
  \title{Independence by random scaling}
  \begin{aug}
    \author{\fnms{Lancelot F.\ }\snm{James}\ead[label=e1]{lancelot@ust.hk}}
    \and
    \author{\fnms{Peter\ }\snm{Orbanz}\corref{}\ead[label=e2]{porbanz@stat.columbia.edu}}
    \affiliation{HKUST and Columbia University}
    \address{Deparment of Information Systems, Business\\ Statistics, and Operations Management\\
      Clear Water Bay, Kowloon\\
      Hong Kong\\
      \printead{e1}
    }
    \address{Department of Statistics\\
      1255 Amsterdam Avenue\\
      New York, NY-10027, USA\\
      \printead{e2}
    }    
  \end{aug}
  \maketitle
\end{frontmatter}

\def\N{M_0}
\def\T{T_{0}}
\def\Tmod{T}
\def\Nmod{M}
\def\ind#1{\text{\tiny #1}}
\def\fT{f_{T_{\!{_0}}}}
\def\fTmod{f_{\ind{\rm T}^{\ast}}}
\def\P{\mathbb{P}}
\def\PK{\text{\rm PK}}
\def\Levy{L\'evy\ }
\def\PD{\text{\rm PD}}
\def\etilt{^\ind{\rm(s)}}
\def\ptilt{^\ind{\rm[$\nu$]}}
\def\condind{{\perp\!\!\!\perp}}
\def\ie{i.e.\ }
\def\eg{e.g.\ }
\def\mean{\mathbb{E}}
\def\equdist{\stackrel{\text{\rm\tiny d}}{=}}
\def\equas{=_{\text{\rm\tiny a.s.}}}
\def\braces#1{{\lbrace #1 \rbrace}}
\def\simiid{\sim_{\mbox{\tiny iid}}}
\def\Law{\mathcal{L}}
\def\iid{i.i.d.\ }

\section{Introduction and main results}

Distributional identities involving elementary
random variables play an important role in probability and related fields.
Such identities arise, for instance, in the study of path
properties of stochastic processes \citep{Pakes:Kathree:1992:1,Matsumoto:Yor:2001:1,
Letac:Wesolowski:2000:1}, and in applications of Stein's method
\citep{Pekoz:Rollin:Ross:2016:1}. A fundamental example is the following:
For ${a>0}$, generically denote by $G_a$ a $\text{Gamma}(a,1)$ variable. If $G_a$ and $G_b$ are
independent, then ${(G_a+G_b)\,\condind\, G_a/(G_a+G_b)}$. \citet{Lukacs:1955} has shown
this property is exclusive to gamma variables, and hence characterizes the gamma distribution.
This result and its ramifications
are collectively known as the \emph{beta-gamma algebra}. Its relevance to path properties of 
Brownian motion and related phenomena is highlighted by \citet{Revuz:Yor:1999:1}.

\def\Gat{G_{\!\frac{\theta}{\raisebox{1pt}{\tiny $\alpha$}}}}
The distributional properties studied in the following are of the form
\begin{equation}
  \label{eq:TXT}
  T\,\condind\, \xi T\qquad\text{ for positive random variables }\xi\text{ and }T\;.
\end{equation}
Lukacs' characterization shows such variables exist---take ${T=G_a+G_b}$ and ${\xi=1/G_a}$---but 
also implies $T$ is a sum of independent variables only if these variables are gamma. 
\citet{Pitman:Yor:1997} have identified another case:
Fix ${\alpha\in(0,1)}$ and ${\theta>0}$, and abbreviate ${\zeta:=G_{\theta/\alpha}}$.
Let $f_{\alpha}$ be an $\alpha$-stable density, $S_{\alpha,\theta}$ a variable
with density proportional to ${t^{-\theta}f_{\alpha}}$, and denote by ${(\tau_{\alpha}(y))_{y\geq 0}}$ a generalized
gamma subordinator, \ie a non-decreasing \Levy process on ${(0,\infty)}$ with \Levy density
${t\mapsto\alpha t^{-\alpha-1}e^{-t}/\Gamma(1-\alpha)}$. Then
\begin{equation}
  \label{eq:stable:gamma:algebra:1}
  \text{(i)}\quad
  \frac{\tau_{\alpha}(\zeta)}{\zeta^{1/\alpha}}
  \,\condind\,
  \tau_{\alpha}(\zeta)
  \qquad\quad\text{(ii)}\quad
  \tau_{\alpha}(\zeta)\equdist G_{\theta} 
  \qquad\quad\text{(iii)}\quad
  \frac{\tau_{\alpha}(\zeta)}{\zeta^{1/\alpha}}\equdist S_{\alpha,\theta}\;,
\end{equation}
which follows from the proof of \citep[][Proposition 21]{Pitman:Yor:1997}. Rescaling
to ${\tilde{\tau}_\alpha(y):=\tau_{\alpha}(y)/y^{1/\alpha}}$ gives
\begin{equation*}
    \text{(i)}\;
    \tilde{\tau}_{\alpha}(\zeta)
    \,\condind\,
    \tilde{\tau}_{\alpha}(\zeta)\cdot\zeta^{1/\alpha}
    \qquad\text{(ii)}\;
    \tilde{\tau}_{\alpha}(\zeta)\cdot \zeta^{1/\alpha}
    \equdist G_{\theta}
    \qquad\text{(iii)}\;
    \tilde{\tau}_{\alpha}(\zeta)
    \equdist
    S_{\alpha,\theta}\;.
\end{equation*}
Clearly, (\ref{eq:stable:gamma:algebra:1}i) is an instance of \eqref{eq:TXT}.

Since both gamma and stable variables are distinguished by their scaling behavior, it 
is natural to ask in how far scaling properties are intrinsic to \eqref{eq:TXT}.
Our first result
shows that the relevant property 
is not scaling per se, but rather a form of exponential tilting. In the case of the stable,
this exponential tilt manifests as a scaling operation; at close inspection, 
the relationship is visible already in \citep{Pitman:Yor:1992:1}.
Let $\T$ be a non-negative random variable with density $\fT$ and cumulant function
${\psi(s):=-\log\mean[e^{-s\T}]}$.
For any such random variable, the exponentially tilted variable $\T\etilt$ and the
polynomially tilted variable $\T\ptilt$ are given by the densities
\begin{equation*}
  \P(\T\etilt\in dt)=e^{-st+\psi(s)}\fT(t)dt
  \qquad\text{ and }\qquad
  \P(\T\ptilt\in dt)=\frac{t^{-\nu}}{\mean[\T^{-\nu}]}\fT(t)dt\;,
\end{equation*}
for ${s\geq 0}$, and for ${\nu>0}$ chosen such that ${\mean[\T^{-\nu}]<\infty}$. 
\begin{theorem}
  \label{theorem}
  Fix ${\nu>0}$, and let $\T$ be a positive random variable with cumulant function
  $\psi$ and ${\mean[\T^{-\nu}]<\infty}$.
  Let $\xi$ and $\Tmod$ be positive, absolutely continuous random variables. Then
  \begin{equation}
    \label{eq:main:2}
    \text{(i)}\quad
    \Tmod\condind\,\xi\Tmod 
    \qquad\quad\text{(ii)}\quad
    \xi\Tmod\equdist G_{\nu} 
    \qquad\quad\text{(iii)}\quad
    {\Tmod\equdist\T\ptilt}
  \end{equation}
  holds if and only if the pair $(\Tmod,\xi)$ satisfies
  \begin{equation}
    \label{eq:main:1}
    \P_{\psi,\nu}(\xi\in ds)
    =
    \frac{e^{-\psi(s)}s^{\nu-1}}{\mean[\T^{-\nu}]\Gamma(\nu)}ds
    \qquad\text{ and }\qquad
    \Tmod|(\xi=s)\equdist\T\etilt \;.
  \end{equation}
\end{theorem}
\noindent Conditional tilting thus yields a large class of
random variables satisfying \eqref{eq:TXT}, and the scaled variable is
always gamma.

The variables in \eqref{eq:stable:gamma:algebra:1} take scalar values.
It is shown in \citep{Pitman:Yor:1997}, however,
that the property extends to the entire path of the process $\tau_{\alpha}$:
For any ${y\in[0,1]}$,
\begin{equation}
  \label{eq:stable:gamma:algebra:2}
  \text{(i)}\quad
  \frac{\tau_{\alpha}(y\zeta)}{\zeta^{1/\alpha}}
  \,\condind\,
  \tau(\zeta)
  \qquad\quad\text{(ii)}\quad
  \tau_{\alpha}(\zeta)\equdist G_{\theta} 
  \qquad\quad\text{(iii)}\quad
  \frac{\tau_{\alpha}(\zeta)}{\zeta^{1/\alpha}}\equdist\tau_{\alpha}(1)\ptilt\;.
\end{equation}
Combined with \cref{theorem}, this suggests an analogous result for general subordinators, which we
state in terms of random measures:
Let $\Omega$ be a Polish space, $\mu$ a probability measure on $\Omega$, and $\lambda$ a \Levy
density on ${(0,\infty)}$. We assume $\lambda$ is strictly positive and continuous, and 
${\int_{0}^{\infty}\min\braces{1,t}\lambda(t)dt<\infty}$.
Let ${(J_n,\omega_n)}$ be the points of a Poisson process on ${(0,\infty)\times\Omega}$ 
with mean measure ${\lambda(t)dt\mu(d\omega)}$. Then
${N:=\sum_{n}J_n\delta_{\omega_n}}$ is a random measure 
on $\Omega$, with ${N(\Omega)<\infty}$ a.s.
If $h$ is a non-negative function with ${\mean[h(N(\Omega))]=1}$,
the random measure $M$ specified by 
\begin{equation}
  \label{eq:weighted:PRM}
  \P(M\in dm)=h(m(\Omega))\P(N\in dm)
\end{equation}
again satisfies ${M(\Omega)<\infty}$. 
If $\psi$ is the cumulant function of the scalar variable ${M(\Omega)}$, 
one can define an exponential tilt ${M\etilt}$
of $M$ as ${\P(M\etilt\!\in dm)=e^{\psi(s)-sm(\Omega)}\P(M\in dm)}$.
\begin{theorem}
  \label{theorem:2}
  Let $\N$ and $\Nmod$ be random measures of the general form \eqref{eq:weighted:PRM}, let
  $\psi$ be the cumulant function of ${\N(\Omega)}$, and fix ${\nu>0}$
  such that ${\mean[\N(\Omega)^{-\nu}]<\infty}$. Then
  \begin{equation}
    \label{eq:theorem:2}
    \text{(i)}\quad
    \Nmod\,\condind\,\xi\Nmod(\Omega)
    \qquad\quad\text{(ii)}\quad
    \xi \Nmod(\Omega)\equdist G_{\nu}
    \qquad\quad\text{(iii)}\quad
    \Nmod(\Omega)\equdist\N(\Omega)^{[\nu]}
  \end{equation}  
  if and only if ${\xi\sim\P_{\psi,\nu}}$ and 
  ${\Nmod|(\xi=s)=\N\etilt}$.
\end{theorem}
\noindent The results are related through the total masses of the random measures: If $\N$ and $\Nmod$
satisfy \cref{theorem:2}, their total masses ${\T:=\N(\Omega)}$ and ${\Tmod:=\Nmod(\Omega)}$
satisfy \cref{theorem}.

\cref{theorem:2} can be applied to normalized random  measures:
As ${T=M(\Omega)}$ is almost surely finite,
${P:=M/T}$ is a random discrete probability measure
\citep{Kingman:1975}.
Since $M$ in 
(\ref{eq:theorem:2}i) is independent of ${\xi T}$, and $T$ is a functional of $M$,
it follows that
\begin{equation}
  \label{eq:independence:rpm}
  P\,\condind\,\xi T
  \qquad\text{ where }\qquad 
  \xi T\equdist G_{\nu}\;.
\end{equation}
The conditions above imply $P$ is of the form 
\begin{equation*}
  P\equdist \sum_{n\in\mathbb{N}}P_n\delta_{\omega_n}
  \qquad\text{ where }\qquad
  (P_n)\condind (\omega_n)
  \quad\text{ and }\quad
  \omega_1,\omega_2\simiid\mu\;,
\end{equation*}
and $(P_n)$ is a random sequence ${P_1\geq P_2\geq\ldots}$ with ${\sum_n P_n=1}$ almost surely.
It is hence no loss of generality to assume $\mu$ is the uniform law on $[0,1]$, or 
to neglect the atoms $\omega_n$ altogether. Throughout, we treat
random probability measures and random sequences $(P_n)$ interchangeably.
If the total mass ${N(\Omega)}$ of the random measure in \eqref{eq:weighted:PRM}
has density $f$, then 
${T=M(\Omega)}$ has density ${hf}$. 
This density, and the \Levy density $\lambda$ of 
$N$, completely determine the law of ${(P_n)}$, which is called a Poisson-Kingman distribution \citep{Pitman:2003:1},
and denoted ${\PK(\lambda,hf)}$.
A distinguished example within the Poisson-Kingman family are the two-parameter Poisson-Dirichlet
distributions ${\PD(\alpha,\theta)}$, with parameters ${\alpha\in[0,1)}$ and ${\theta>-\alpha}$
\citep{Kingman:1975,Pitman:Yor:1997}. 
For ${\theta>0}$, they are known to satisfy \eqref{eq:independence:rpm}: 
There is a random measure $\Nmod$ with total mass $\Tmod$ such that ${P=\Nmod/\Tmod}$ satisfies
\begin{equation*}
  P\,\condind\, \xi\Tmod
  \quad\text{ and }\quad
  \xi\Tmod\equdist G_{\theta+\alpha}
  \quad\text{ if }\quad
  (P_n)\sim\PD(\alpha,\theta) \;.
\end{equation*}
If ${\alpha>0}$, this is once again Proposition 21 of 
\citep{Pitman:Yor:1997}, and can be derived from \eqref{eq:stable:gamma:algebra:2} by choosing
${M[0,y]:=\tau_{\alpha}(y\zeta)/\zeta^{1/\alpha}}$, in which case ${(P_n)}$ has law 
${\PD(\alpha,\theta)}$. 
If ${\alpha=0}$, choose ${M[0,y]=\tau(y\theta)}$ for a gamma subordinator $\tau$ instead;
then ${(P_n)\sim\PD(0,\theta)}$, and the result follows from Lukacs' characterization.

Relative to Proposition 21 of \citet{Pitman:Yor:1997},
our results imply an extension to the case ${\theta\leq 0}$:
Start with a generalized gamma subordinator ${\tau_{\alpha}(y)}$ and ${b>0}$. Size-biasing the 
process ${\tau_{\alpha}(yb^{\alpha})/b}$ turns it into a bridge
\begin{equation*}
  \frac{\tau_{\alpha}(yb^{\alpha})}{b}+\frac{G_{1-\alpha}}{b}\mathds{1}\braces{U\leq y}
  \qquad\text{ for }U\sim\text{Uniform}[0,1]\text{ independently.}
\end{equation*}
In \cref{sec:PD}, we construct scalar random variables ${H_{\alpha,\theta}}$ and 
${\xi_{H_{\alpha,\theta}}}$ such that randomizing
$b$ by ${\xi_{H_{\alpha,\theta}}+H_{\alpha,\theta}}$ defines a random measure
\begin{equation}
  \label{eq:PY:bridge}
  M[0,y]:=\frac{\tau_{\alpha}(y(\xi_{H_{\alpha,\theta}}+H_{\alpha,\theta})^{\alpha})}{\xi_{H_{\alpha,\theta}}+H_{\alpha,\theta}}
  +\frac{G_{1-\alpha}}{\xi_{H_{\alpha,\theta}}+H_{\alpha,\theta}}\mathds{1}\braces{U\leq y}
\end{equation}
for which the weights $(P_n)$ of ${P=\Nmod/\Tmod}$ have law ${\PD(\alpha,\theta)}$.
\begin{proposition}
  \label{result:PY:full:range}
  For any ${\alpha\in[0,1)}$ and ${\theta>-\alpha}$, the ranked weights $(P_n)$ derived from \eqref{eq:PY:bridge}
    satisfy  ${P\sim\PD(\alpha,\theta)}$, independently of 
    \begin{equation*}
      (\xi_{H_{\alpha,\theta}}+H_{\alpha,\theta}) \Tmod\equdist G_{1-\theta}
      \qquad\text{ and of }\qquad
    \xi_{H_{\alpha,\theta}} \Tmod\equdist G_{\theta+\alpha}\;.
    \end{equation*}
\end{proposition}
The remainder of this article describes applications and examples; proofs are collected in the appendix. 
While the results apply to quite general processes, our examples emphasize the stable subordinator,
which leads to interesting extensions of Proposition 21 of \citep{Pitman:Yor:1997}.

\section{Application to generalized gamma subordinators}
\label{sec:GG}

In this section, we consider the scaled and time-changed 
process ${\tau_{\alpha}(yb^{\alpha})/b}$ that already arose above.
This process can be equivalently represented by exponentially tilting a stable subordinator
\citep{Pitman:Yor:1997}:
Let $f_{\alpha}$ denote the density of an $\alpha$-stable random variable, and ${\lambda_{\alpha}}$ 
an $\alpha$-stable \Levy density. If $\sigma^{(b)}_{\alpha}$ is a subordinator with \Levy density 
${e^{-bt}\lambda_{\alpha}(t)}$, then
\begin{equation*}
  \sigma^{(b)}_{\alpha}(y)\equdist \frac{\tau_{\alpha}(yb^{\alpha})}{b}\;,
\end{equation*}
where the left-hand side is well-defined even if ${b=0}$.
The variable
\begin{equation*}
  X_{\alpha,b}:=\frac{\tau_{\alpha}(b^{\alpha})}{b}
  \quad\text{ hence has density }\quad
  t\mapsto e^{-bt+b^\alpha}f_{\alpha}(t)\;.
\end{equation*}

Exponentially tilted \Levy densities as the one above define a class of Poisson-Kingman distributions 
for which our results take a special form:
For a \Levy density $\lambda$ and ${\nu>0}$, let $\T$ be the total mass of a random measure
defined by $\lambda$. The Poisson-Kingman distribution ${\PK(\lambda,\Law(\T\ptilt))}$
can be embedded in a one-parameter family
${\PK(e^{-bt}\lambda(t),\Law(T_b\ptilt))}$, for ${b\geq 0}$, where $T_b$ is 
the total mass of a random measure defined by ${e^{-bt}\lambda(t)}$.
The conditioning operation in \cref{theorem:2} then takes the form of a parameter shift:
A random probability measure ${P}$ with law ${\PK(e^{-bt}\lambda(t),\Law(T_b^{[\nu]}))}$ 
satisfies \eqref{eq:independence:rpm} 
if and only if 
\begin{equation*}
  {\Tmod|(\xi=s)\equdist\Tmod_{b+s}}
  \quad\text{ and }\quad
            {\xi\sim\P_{\psi,\nu}}
\end{equation*}
where $\psi$ is the cumulant function of $T_{b=0}$.

\subsection{The basic case}  
\label{sec:gg}

Suppose the random measure $\N$ in \cref{theorem:2} is defined as
${\N[0,y]:=\tau_{\alpha}(yb^\alpha)/b}$ for all ${y\in[0,1]}$.
The total mass ${\N[0,1]\equdist X_{\alpha,b}}$
then has cumulant function ${\psi(s)=(b+s)^{\alpha}-b^{\alpha}}$, and substituting into
the theorem yields
\begin{equation*}
  \P_{\psi,\nu}(\xi\in ds)=
  \frac{e^{-(b+s)^{\alpha}+b^{\alpha}}s^{\nu-1}}{\mean[X_{\alpha,b}^{-\nu}]\Gamma(\nu)}ds
  \qquad\text{ and }\qquad
  \Nmod[0,y]\big\vert(\xi=s)
  \;
  =
  \;
  \frac{\tau_{\alpha}(y(b+s)^\alpha)}{b+s}\;.
\end{equation*}
Consequently, the random probability measure
\begin{equation}
  \label{eq:example:gg:independence}
  P[0,y]:=\frac{\tau_{\alpha}(y(b+\xi)^\alpha)}{\tau_{\alpha}((b+\xi)^\alpha)}
  \qquad\text{ satisfies }\qquad
  P\;\condind \;\xi\frac{\tau_{\alpha}((b+\xi)^\alpha)}{b+\xi}\;
  \equdist\;
  G_{\nu}\;.
\end{equation}
The variables ${\Tmod=\Nmod[0,1]}$ and ${\T=\N[0,1]}$ satisfy
\begin{equation}
  \label{eq:total:mass:stable}
  \Tmod=\T\ptilt\equdist\frac{\tau_{\alpha}((b+\xi)^{\alpha})}{b+\xi}
  \qquad\text{ and hence }\qquad
  \mathbb{P}(T\in dt)\propto t^{-\nu}e^{-bt+b^{\alpha}}f_{\alpha}(t)dt\;.
\end{equation}
The resulting law of ${(P_n)}$ is ${\PK(\lambda_{\alpha},\Law(\Tmod))}$.
For ${b:=0}$ and ${\nu>0}$,
this law is specifically $\PD(\alpha,\nu)$, which recovers Proposition 21 of 
\citet{Pitman:Yor:1997}.
Both the independence property in \eqref{eq:example:gg:independence} \emph{and}
equality in distribution to $G_\nu$ remain true if 
$b$ is randomized by mixing against any positive random variable.

\subsection{Size-biasing}
\label{sec:sb:gg}

If $Y$ is any positive random variable with density $f_{\ind{Y}}$, 
we denote by $Y^{\ast}$ the size-biased variable with density
${yf_{\ind{Y}}(y)/\mean[Y]}$. 
For an independent uniform variable $U$ on ${[0,1]}$, the process
\begin{equation*}
  \tau_{\alpha,b}^{\ast}(y):=\frac{\tau_{\alpha}(yb^{\alpha})}{b}+\frac{G_{1-\alpha}}{b}\mathds{1}\braces{U\leq y}
  \qquad\text{ hence }\qquad
  \tau_{\alpha,b}^{\ast}(1)
  \equdist
  \Bigl(\frac{\tau_{\alpha}(b^{\alpha})}{b}\Bigr)^{\ast}\;,
\end{equation*}
and can be regarded as a size-biased form of ${\tau_\alpha(yb^\alpha)/b}$
\citep{Pakes:Sapatinas:Fosam:1996:1,Perman:Pitman:Yor:1992}.
Since the summands are independent, their cumulant functions are additive, and 
the cumulant function of ${\tau^{\ast}_{\alpha,b}(1)}$ is 
\begin{equation*}
\psi(s)=-(\alpha-1)\log(1+{\textstyle\frac{s}{b}})+(b+s)^{\alpha}-b^{\alpha}\;.
\end{equation*}
For the random measure defined on the interval by ${\N[0,y]:=\tau^{\ast}_{\alpha,b}(y)}$,
the distributions in \cref{theorem:2} then take the form 
\begin{equation}
  \label{eq:xi:sbgg}
  \P_{\psi,\nu}(\xi\in ds)=
  \frac{\alpha(b+s)^{\alpha-1}e^{-(b+s)^{\alpha}+b^{\alpha}}s^{\nu-1}}{\Gamma(\nu)\mean[\tau_{\alpha,b}(1)^{-\nu+1}]}ds
  \qquad\text{ and }\qquad
  \Nmod\big\vert(\xi=s)=\tau_{\alpha,b+s}^{\ast}\;.
\end{equation}
As the variable ${\tau_{\alpha,b}^{\ast}(1)}$ can be defined by tilting and size-biasing a stable
variable, its density is
${g_{\alpha,b}(t):=b^{1-\alpha}e^{-bt+b^{\alpha}}tf_{\alpha}(t)/\alpha}$. The marginal law of ${T=M[0,1]}$ is then
\begin{equation*}
  \P(T\in dt)
  =\frac{t^{-\nu}g_{\alpha,b}(t)}{\int s^{-\nu}g_{\alpha,b}(s)ds}dt
  =\P(X_{\alpha,b}^{[\nu-1]}\in dt)\;,
\end{equation*}
and we obtain:
\begin{proposition}
  \label{prop:size:bias}
  Let $f_{\alpha}$ be the $\alpha$-stable density.
  If the weights of a random probability measure $P$ have law ${\PK(f_\alpha,\Law(X_{\alpha,b}\ptilt))}$, it
  can be represented as ${P=M/T}$ for
  \begin{equation*}
    \Nmod[0,y]
    =
    \frac{\tau_{\alpha}(y(b+\xi)^{\alpha})}{b+\xi}+\frac{G_{1-\alpha}}{b+\xi}\mathds{1}\braces{U\leq y}
    =
    \frac{\tau_{\alpha}(y((b+\xi)^{\alpha}+G_{(1-\alpha)/\alpha}\mathds{1}\braces{U\leq y}))}{b+\xi}
\end{equation*}
  and satisfies
  ${P\condind\,\xi\Tmod}$ and ${\xi\Tmod\equdist G_{\nu}}$
  and ${\Tmod\equdist X_{\alpha,b}^{[\nu-1]}}$.
\end{proposition}
\noindent For example, choose ${\nu=1}$, and abbreviate ${Z:=(G_1+b^{\alpha})^{1/\alpha}}$. Then, for any ${b>0}$,
\begin{equation*}
  \frac{\tau_{\alpha}(Z^{\alpha}+G_{\frac{1-\alpha}{\alpha}})}{Z}
  \quad\condind\quad
  (Z-b)
  \frac{\tau_{\alpha}(Z^{\alpha}+G_{\frac{1-\alpha}{\alpha}})}
  {Z}
  \quad
  \equdist
  \quad
  G_1\;,
\end{equation*}
\ie the value of the process $\tau_{\alpha}$, taken at a suitable random time,
decouples from itself by random scaling. This also shows the unbiased case
in \cref{sec:gg} can be recovered from the size-biased one by choosing ${\nu=1}$:
Observe the term on the left is distributed as
${\tau_{\alpha}(Z^{\alpha}+G_{\frac{1-\alpha}{\alpha}})/Z\equdist\tau_{\alpha}(b^{\alpha})/b}$,
for any ${b>0}$. 
For ${b'\geq 0}$ and ${\nu'>0}$, we may substitute
${b=b'+\xi'}$, where ${\xi'}$ has
density proportional to ${e^{-(b'+s)^{\alpha}+b'^{\alpha}}s^{\nu'-1}}$ as in \cref{sec:gg} above.
Then
\begin{equation*}
  \Tmod\equdist
  \frac{\tau_{\alpha}((b'+\xi')^{\alpha})}{(b'+\xi')}
  \qquad\text{ and hence }\qquad
  \mathbb{P}(\Tmod\in dt)\propto t^{-\nu'}e^{-b't+b'^{\alpha}}\;,
\end{equation*}
which recovers all cases in \cref{sec:gg}.

\subsection{Poisson-Dirichlet models}
\label{sec:PD}

A $\PD(\alpha,\theta)$ random measure $P$ can be represented as
\begin{equation}
  \label{eq:Dong:et:al}
  P[0,y]
  \equdist
  \frac{\tau_{\alpha}(G_{(\alpha+\theta)/\alpha}y+G_{(1-\alpha)/\alpha}\mathbb{I}\braces{U\leq y})}
        {\tau_{\alpha}(G_{(\alpha+\theta)/\alpha}+G_{(1-\alpha)/\alpha})}
        \qquad\text{ for any }\theta>-\alpha\;.
\end{equation}
This can be read from \citet*{Dong:Goldschmidt:Martin:2006:1}, 
or indeed from \citet{Pitman:Yor:1992:1}.
Define 
\begin{equation*}
  H_{\alpha,\theta}:=G_{\frac{\theta+\alpha}{\alpha}}^{1/\alpha}B_{1-\alpha,\theta+\alpha}\;.
\end{equation*}
Now index the random variable $\xi$ in \eqref{eq:xi:sbgg} explicitly by the value of $b$ as $\xi_b$,
and let ${\xi_{H_{\alpha,\theta}}}$ denote the variable obtained by mixing $b$ against ${H_{\alpha,\theta}}$.
\begin{lemma}
  \label{lemma:beta:gamma}
  Let $\P_{\psi,\nu}$ be defined as in \eqref{eq:xi:sbgg}. For each ${b>0}$, let 
  ${\xi_b\sim\P_{\psi,\theta+\alpha}}$. Then
  \begin{equation}
    \label{eq:beta:gamma}
    \begin{split}
    \text{(i)}\quad&
    (\xi_{H_{\alpha,\theta}}+H_{\alpha,\theta})^{\alpha}\equdist G_{(\alpha+\theta)/\alpha}\\
    \text{(ii)}\quad&
    \bigl(
    \xi_{H_{\alpha,\theta}},H_{\alpha,\theta}
    \bigr)
    \equdist
    G_{\frac{\theta+\alpha}{\alpha}}^{1/\alpha}(B_{\theta+\alpha,1-\alpha},1-B_{\theta+\alpha,1-\alpha})\;.
    \hspace{2cm}
    \end{split}
  \end{equation}
\end{lemma}
\noindent We have hence established the result stated in the introduction:
\begin{proof}[Proof of \cref{result:PY:full:range}]
  Substitution of ${\xi_{H_{\alpha,\theta}}+H_{\alpha,\theta}}$ for $b$ in \cref{prop:size:bias}
  yields the random measure defined in \eqref{eq:PY:bridge}. By \eqref{eq:Dong:et:al}, it 
  normalizes to a measure with weights ${(P_n)\sim\PD(\alpha,\theta)}$, and 
  the claim follows from \cref{prop:size:bias}.
\end{proof}

\subsection{Implications for $\alpha$-diversities}
For ${(P_n)\sim\PK(\lambda_{\alpha},hf_{\alpha})}$, it is known \citep{Pitman:2006} that 
\begin{equation*}
  \Gamma(1-\alpha)^{-1}\lim_{\varepsilon\rightarrow\infty}\varepsilon^{\alpha}|\braces{n|P_n\geq\varepsilon}|=T^{-\alpha}
  \qquad\text{a.s.}
\end{equation*}
The random variable $T^{-\alpha}$ can be interpreted in terms of a local time, and is also
known as the $\alpha$-diversity of the exchangeable random partition of $\mathbb{N}$ defined by $P$: 
If $K_n$ is the number of distinct
blocks in the restriction of this partition to the subset $[n]$, then
${n^{-\alpha}K_n\rightarrow T^{-\alpha}}$ almost surely as ${n\rightarrow\infty}$.
The case ${T=S_{\alpha,\theta}}$ arises in Bayesian statistics,
stochastic processes, and models for random trees and graphs 
\citep[e.g.][]{Favaro:Lijoi:Mena:Pruenster:2009,Goldschmidt:Haas:2015,James:2007uq,Pekoz:Rollin:Ross:2016:1,Winkel:2005:1}.

The law considered above is ${(P_n)\sim\PK(\lambda_{\alpha},\Law(T))}$, 
where ${T=\tau_{\alpha}((b+\xi)^{\alpha})/(b+\xi)}$ as
in \eqref{eq:total:mass:stable}. 
For ${b=0}$, the variable $T^{-\alpha}$ is the $\alpha$-diversity of the two-parameter Chinese restaurant process
\citep{Pitman:2006}. More generally, for any value ${b\geq 0}$, the resulting partition is 
of Gibbs type \citep{Pitman:2006},
since $\lambda_{\alpha}$ defines a stable subordinator. There hence exists a subclass of Gibbs-type
measures that is strictly larger than the Poisson-Dirichlet family, and whose $\alpha$-diversity
exhibits a similar independence property ${\xi^{-\alpha}T^{-\alpha}\condind\, P}$.

\section{Application to excursions straddling an exponential random time}
\label{sec:BFRY}

This section considers applications to a type of distributions and processes
that arise in a range of contexts, including
passage times of \Levy processes, excursions of regular linear diffusions, interval partitions
generated by a subordinator, and also in applications in statistics and finance
\citep[e.g.][]{Bertoin:Fujita:Roynette:Yor:2006:1,Gnedin:Pitman:2005:1,James:Orbanz:Teh:2015:1, 
Kyprianou:2006:1,Pitman:Yor:1997,Pitman:1997:1,Winkel:2005:1}.

\subsection{Independence of scaled excursion durations}
Let $\tau$ again be a subordinator, with \Levy density $\lambda$ and ${\tau(1)<\infty}$
a.s., and denote its \Levy exponent
${\Psi(s):=\int_0^{\infty}(1-e^{-st})\lambda(t)dt}$.
Following \citet{Winkel:2005:1}
and the exposition in \citep{Bertoin:Fujita:Roynette:Yor:2006:1},
define the local time process $L$, overshoot process $O$, and undershoot
process $U$ as
\begin{equation*}
  L_t:=\inf\lbrace s\,\vert\,\tau(s)>t\rbrace
  \qquad\quad
  O(t,\tau):=\tau({L_{t}})-t
  \qquad\quad
  U(t,\tau):=t-\tau({L_{t}}-)\;,
\end{equation*}
where ${\tau({L_{t}}-)}$ is the prepassage height, \ie the left-hand limit ${\lim_{t\nearrow L_t}\tau(t)}$.
For an independent exponential time $G_1$, abbreviate
\begin{equation*}
  O(\tau):=O(G_1,\tau)
  \quad\text{ and }\quad
  U(\tau):=U(G_1,\tau)
  \qquad\text{ and define }\qquad
  \Delta(\tau):=  O(\tau)+U(\tau)\;.
\end{equation*}
The variable ${\Delta(\tau)}$
can be interpreted as the duration of the excursion from $0$ to $0$ of a strongly recurrent linear 
diffusion that straddles the random time $G_1$, and whose inverse local time is $\tau$
\cite{Salminen:Vallois:Yor:2007:1}.
The density $f_\lambda$ of $\Delta(\tau)$ and the joint density $h_{\lambda}$ of ${(O(\tau),U(\tau))}$ are known to be
\begin{equation*}
  f_{\lambda}(t)
  :=
  \frac{(1-e^{-t})\lambda(t)}{\Psi(1)}
  \qquad\text{ and }\qquad
  h_{\lambda}(v,w):=\frac{e^{-v}\lambda(v+w)}{\Psi(1)}\;,
  \label{straddle,straddlejoint}
\end{equation*}
see \cite{Pitman:1997:1,Winkel:2005:1,Salminen:Vallois:Yor:2007:1}
regarding $f_{\lambda}$, and \cite[][eq. (33)]{Salminen:Vallois:Yor:2007:1} for $h_\lambda$.
For ${b\geq 0}$, we define $\tau^{(b)}$ as the subordinator with exponentially tilted \Levy density
\begin{equation*}
  \lambda^{(b)}(s)=e^{-bs}\lambda(s)
  \qquad\text{ which has \Levy exponent }\qquad
  \Psi(s+b)-\Psi(b)\;.
\end{equation*}
An additional polynomial tilt yields the subordinator $\tau^{(b)}_{\nu}$ with \Levy density
\begin{equation*}
  \lambda^{(b)}_{\nu}(s)=s^{\nu}e^{-bs}\lambda(s)
  \quad\text{ with \Levy exponent }\quad
  \Psi_{\nu}(b):=\int_{0}^{\infty}(1-e^{-s})\lambda^{(b)}_{\nu}(s)ds\;.
\end{equation*}
If the scalar variable $\T$ in \cref{theorem} is chosen as ${\T:=\Delta(\tau^{(b)}_{\nu})}$, the resulting
law of $\xi$ is
\begin{equation}
  \label{eq:local:time:xi}
  \P_{\psi,\nu}(\xi\in ds)=\frac{\Psi_{\nu}(b+s)s^{\nu-1}}{\Gamma(\nu)(\Psi(b+1)-\Psi(b))}ds\;.
\end{equation} 
The variable $\Tmod$ in the theorem is then ${\Tmod=\Delta(\tau^{(b+\xi)}_{\nu})}$.
\begin{proposition} 
  \label{proposition:BFRY}
  Fix ${\nu>0}$ and ${b\geq 0}$, let $\xi$ be a random variable with law \eqref{eq:local:time:xi}.
  Then the conditional density of ${\Law(O(\tau_{\nu}^{(b+\xi)}),U(\tau_{\nu}^{(b+\xi)})|\xi=s)}$ is
  $h_{\lambda_{\nu,b+s}}$, and
  \begin{equation*}
    (O(\tau^{(b+\xi)}_{\nu}),U(\tau^{(b+\xi)}_{\nu}))
    \equdist
    (O(\tau^{(b)}),U(\tau^{(b)}))
    \qquad\text{ independently of }\qquad
    \xi\Delta(\tau^{(b+\xi)}_{\nu})\equdist G_{\nu}\;.
  \end{equation*}
  The process ${\tau^{(b+\xi)}_{\nu}|\xi=s}$ is compound Poisson with rate 
  and jump density
  \begin{equation*}
    r_{b,\nu}:=\int_{0}^{\infty}{\mbox e}^{-(b+s)t}t^{\nu}\lambda(t)dt
    \quad\text{ and }\quad
    s\mapsto 
    e^{-(b+s)t}t^{\nu}\lambda(t)/r_{b,\nu}
  \end{equation*}
  whenever $r_{b,\nu}<\infty$, in particular for ${\nu \ge 1}$.
\end{proposition}
\noindent Since ${\Delta=O+U}$, it follows that the excursion duration satisfies
\begin{equation*}
  \Delta(\tau^{(b+\xi)}_{\nu})\equdist\Delta(\tau^{(b)})
  \qquad\text{ independently of }\qquad
  \xi\Delta(\tau^{(b+\xi)}_{\nu})\;.
\end{equation*}
The result does not imply independence of $\xi\Delta(\tau^{(b+\xi)}_{\nu})$ and the entire process ${\tau^{(b+\xi)}_{\nu}}$.

\subsection{A concrete example}
Let $\tau^{(b+\xi)}_{\nu}$ have \Levy density 
\begin{equation}
  \label{eq:density:BFRY:gg}
  \lambda^{(b+\xi)}_{\nu}(s)=\frac{\alpha}{\Gamma(1-\alpha)}s^{\nu-\alpha-1}e^{-b s}
         \qquad\text{ for }s\in(0,\infty)\;.
\end{equation}
Changing parameters to ${\delta:=\alpha-\nu}$, and comparing to the generalized
gamma subordinator $\tau_{\alpha}$ used in \cref{sec:GG}, shows \eqref{eq:density:BFRY:gg}
is up to a constant the \Levy density of the subordinator ${\tau_{\delta}(tb^{\delta})/b}$.
The \Levy exponent $\Psi_{\nu}$ of
$\lambda^{(b)}_{\nu}$, and hence the variable $\xi$ defined by \eqref{eq:local:time:xi},
depend on the sign of $\delta$. We must distinguish three cases:

\begin{enumerate}
\item ${\delta\in(0,\alpha)}$ and ${b\geq 0}$: $\tau^{(b)}_{\nu}$ is a generalized gamma process with infinite
  activity and parameter $\delta$, with
  \begin{equation*}
    \Psi_{\nu}(b)=\frac{\alpha\Gamma(1-\delta)}{\delta\Gamma(1-\alpha)}\bigl({(b+1)}^{\delta}-{b}^{\delta}\bigr)
    \qquad\text{ and }\qquad
    \tau^{(b+\xi)}_{\nu}(t)\equdist
    \frac{\tau_{\delta}(t{(b+\xi)}^{\delta}\frac{\alpha\Gamma(1-\delta)}{\delta\Gamma(1-\alpha)})}{b+\xi}\;,
  \end{equation*}
  where $\tau_{\delta}=\tau_{\alpha=\delta}$ is a generalized gamma subordinator.
\item ${\delta=0}$ and ${b>0}$: $\tau^{(b)}_{\nu}$ is a gamma process, with
  \begin{equation*}
    \Psi_{\nu}(b)=\frac{\alpha}{\Gamma(1-\alpha)}\log(1+1/b))
    \qquad\text{ and }\qquad
    \tau^{(b+\xi)}_{\nu}(t)\equdist
    \frac{\gamma(t \frac{\alpha}{\Gamma(1-\alpha)}))}{b+\xi}\;,
  \end{equation*}
  where $\gamma$ is a gamma subordinator with \Levy density ${s^{-1}e^{-s}}$.
  The weights of the random measure 
  ${P=\tau^{(b+\xi)}_{\nu}(y)/\tau^{(b+\xi)}_{\nu}(1)}$ have law ${\PD(0,\delta)}$, independently of $\xi$.
\item For ${\delta<0}$ and ${b>0}$, one obtains a compound Poisson process, with
  \begin{equation*}
    \Psi_{\nu}(b)=\frac{\alpha\Gamma(\nu-\alpha)}{\Gamma(1-\alpha)}\bigl({b}^{\alpha-\nu}-{(b+1)}^{\alpha-\nu}\bigr)
    \qquad\text{ and }\qquad
    \tau^{(b+\xi)}_{\nu}(t)\equdist\sum_{i=1}^{\tilde{N}(t)}\frac{G_{i}}{b+\xi}\;,
  \end{equation*}
  where $\tilde{N}$ is a Poisson process with rate $(b+\xi)^{\alpha-\nu}\frac{\alpha\Gamma(\nu-\alpha)}{\Gamma(1-\alpha)}$,
  and the variables ${G_1,G_2,\ldots}$ are \iid ${G_i\equdist G_{\nu-\alpha}}$.
\end{enumerate}
In each case, the excursion duration is conditionally distributed as
\begin{equation*}
  \Law(\Delta(\tau^{(b+\xi)}_{\nu})|\xi=s)
  =
  \frac{\alpha(1-{\mbox e}^{-t}){\mbox e}^{-(b+s)t}t^{\nu-\alpha-1}}{\Gamma(1-\alpha)\Phi_{\nu}(b+s)}dt\;,
\label{straddleBFRYgen}
\end{equation*}
and satisfies ${\Delta(\tau^{(b+\xi)}_{\nu})\equdist\Delta(\tau^{(b)}_{0})}$,
independently of ${\xi \Delta(\tau^{(b+\xi)}_{\nu})}$.

\vspace{1cm}
{\noindent\textbf{Acknowledgments.}
  LFJ was supported in
  part by grant RGC-HKUST 601712 of the \mbox{HKSAR}.
  PO was supported in part by grant FA9550-15-1-0074 of AFOSR.
}

\newpage

\bibliography{refs-local}
\bibliographystyle{abbrvnat}

\appendix

\section*{Proofs}

\begin{proof}[Proof of \cref{theorem}]
  If \eqref{eq:main:1} holds, the joint density of ${(\xi,T)}$ is
  \begin{equation*}
    \P(\xi\in ds, T\in dt)
    =
    \frac{e^{-st}s^{\nu-1}}{\mean[\T^{-\nu}]\Gamma(\nu)}f_0(t) ds dt\;,
  \end{equation*}
  if $f_0$ is the density of $\T$. It can be disintegrated either into ${\Law(T|\xi)}$ and ${\Law(\xi)}$, which
  recovers \eqref{eq:main:1}, or into ${\Law(\xi|T)}$ and ${\Law(T)}$, in which case
  \begin{equation}
    \label{eq:proof:main:1}
    \P(T\in dt)=\frac{t^{-\nu}f_0(t)}{\mean[\T^{-\nu}]}dt=\P(\T\ptilt\in dt)\;,
  \end{equation}
  which is just (\ref{eq:main:2}iii).
  For any measurable functions $g$ and $h$, a change
  of variables ${y:=st}$ then yields
  \begin{equation*}
    \begin{split}
    \mean[g(T^{\ast})h(\xi T^{\ast})]
    =&\ \iint g(t)h(ts)\frac{e^{-st}e^{\psi(s)}\fT(t)e^{-\psi(s)}s^{\nu-1}}{\mean[T^{-\nu}]\Gamma(\nu)}dsdt\\
    =&\ \int g(t)\frac{\fT(t)t^{-\nu}}{\mean[T^{-\nu}]}dt
    \int h(y)\frac{e^{-y}y^{\nu-1}}{\Gamma(\nu)}dy \;,
    \end{split}
  \end{equation*}
  so (\ref{eq:main:2}i) and (\ref{eq:main:2}ii) are also true.
  Conversely, assume \eqref{eq:main:2} holds, and hence in particular \eqref{eq:proof:main:1}.
  The joint density of $(\Tmod,\xi\Tmod)$ is then
  \begin{equation*}
    \frac{t^{-\nu}\fT(t)}{\mean[\T^{-\nu}]}dt
    \frac{y^{\nu-1}e^{-y}}{\Gamma(\nu)}dy
    =
    \frac{t^{-\nu}\fT(t)}{\mean[\T^{-\nu}]}dt
    \frac{(st)^{\nu-1}e^{-st}}{\Gamma(\nu)}tds
    =
    \frac{s^{\nu-1}}{\Gamma(\nu)\mean[\T^{-\nu}]}ds\cdot
    e^{-st}\fT(t)dt\;.
  \end{equation*}
  If additionally $h$ is any positive function and ${\psi(s)=\log(h(s))}$,
  \begin{equation*}
    \frac{s^{\nu-1}}{\Gamma(\nu)\mean[\T^{-\nu}]}\frac{h(s)}{h(s)}ds\cdot
    e^{-st}\fT(t)dt
    =
    \frac{s^{\nu-1}e^{-\psi(s)}}{\Gamma(\nu)\mean[\T^{-\nu}]}ds\cdot
    e^{-(st-\psi(s))}\fT(t)dt\;,
  \end{equation*}
  which is the product of the two terms in \eqref{eq:main:1}.
\end{proof}

\noindent The proof of \cref{theorem:2} is similar.

\begin{proof}[Proof of \cref{lemma:beta:gamma}]
  ${H_{\alpha,\theta}}$ has density
  ${b\mapsto b^{-\alpha}e^{-b^{\alpha}}\mean[X_{\alpha,b}^{-\theta-\alpha+1}]/(\Gamma(1-\alpha)\mean[S_{\alpha}^{-(\theta)}])}$.
  Integrating against the density of $\xi_b$ given in \eqref{eq:xi:sbgg} shows ${\xi_{H_{\alpha,\theta}}}$ has
  marginal density
  \begin{equation*}
    s\mapsto Cs^{\theta+\alpha-1}\int_0^1 e^{-y(b+1)^{\alpha}}(b+1)^{\alpha-1}b^{-\alpha}db\;,
    \qquad\text{ hence }\qquad
    \xi_{H_{\alpha,\theta}}\equdist G_{\frac{\theta+\alpha}{\alpha}}B_{1-\alpha,\theta+\alpha}\;.
  \end{equation*}
  Taking Laplace transforms yields \eqref{eq:beta:gamma}(i), which implies
  ${(\xi_{H_{\alpha,\theta}},\xi_{H_{\alpha,\theta}}+H_{\alpha,\theta})^{\alpha}}$
  is equal in distribution to ${G_{(\alpha+\theta)/\theta}(B_{\theta+\alpha,1-\alpha},1)}$,
  and hence yields \eqref{eq:beta:gamma}(ii).
\end{proof}

\begin{proof}[Proof of \cref{proposition:BFRY}]
  (i) holds by construction and \cref{theorem}. To obtain (ii) and (iii),
  abbreviate ${(O_{\xi},U_{\xi}):=(O(\tau_{\nu, b+\xi}),U(\tau_{b+\xi,\nu}))}$ 
  and ${\Delta_{\xi}:=\Delta(\tau_{b+\xi,\nu})}$. The joint density of
  ${(O_\xi,U_\xi,\xi)}$ is then 
  {\scriptsize
  \begin{equation*}
    \frac{e^{-v}e^{-(b+s)(v+w)}
      {(v+w)}^{\nu}
      \lambda(v+w)}{\Psi_{\nu}(b+s)}\times\frac{\Psi_{\nu}(b+s)s^{\nu-1}}{\Gamma(\nu)(\Psi(b+1)-\Psi(b))}
    \quad=\quad
    \frac{{\mbox e}^{-v}{\mbox e}^{-(b+s)(v+w)}
      {(v+w)}^{\nu}
      \lambda(v+w)s^{\nu-1}}{\Gamma(\nu)(\Psi(b+1)-\Psi(b))}\;.
  \end{equation*}
  }
  It follows that, for any measureable function $h$,
  \begin{equation*}
    \mean[h(O_{\xi},U_{\xi})e^{-\omega \xi \Delta_{\xi}}]
    =
    \frac{h(v,w){\mbox e}^{-v}{\mbox e}^{-b(v+w)}{(v+w)}^{\nu}
      {\mbox e}^{-s(1+\omega)(x+v)}\lambda(v+w)s^{\nu-1}}{\Gamma(\nu)\Psi_{\lambda}(b+1)-\Psi_{\lambda}(b)}dsdvdw\;,
  \end{equation*}
  and integrating out $s$ yields 
  ${\mean[h(O_{\xi},U_{\xi})e^{-\omega \xi \Delta_{\xi}}]=(1+\omega)^{-\nu}\mathbb{E}[h((O(\tau_{b}),U(\tau_{b}))]}$.
\end{proof}

\end{document}